# Constructions of Spherical 3-Designs

Béla Bajnok

Department of Mathematics and Computer Science, Gettysburg College, Gettysburg, PA 17325-1486, USA. e-mail:bbajnok@cc.gettysburg.edu

**Abstract.** Spherical $t$-designs are Chebyshev-type averaging sets on the $d$-sphere $S^d \subset R^{d+1}$ which are exact for polynomials of degree at most $t$. This concept was introduced in 1977 by Delsarte, Goethals, and Seidel, who also found the minimum possible size of such designs, in particular, that the number of points in a 3-design on $S^d$ must be at least $n \geq 2d + 2$. In this paper we give explicit constructions for spherical 3-designs on $S^d$ consisting of $n$ points for $d = 1$ and $n \geq 4$; $d = 2$ and $n = 6, 8, \geq 10$; $d = 3$ and $n = 8, \geq 10$; $d = 4$ and $n = 10, 12, \geq 14$; $d \geq 5$ and $n \geq 5(d+1)/2$ odd or $n \geq 2d + 2$ even. We also provide some evidence that 3-designs of other sizes do not exist. We will introduce and apply a concept from additive number theory generalizing the classical Sidon-sequences. Namely, we study sets of integers $S$ for which the congruence $\varepsilon_1 x_1 + \varepsilon_2 x_2 + \cdots + \varepsilon_t x_t \equiv 0 \mod n$, where $\varepsilon_i = 0, \pm 1$ and $x_i \in S$ ($i = 1, 2, \ldots, t$), only holds in the trivial cases. We call such sets Sidon-type sets of strength $t$, and denote their maximum cardinality by $s(n, t)$. We find a lower bound for $s(n, 3)$, and show how Sidon-type sets of strength 3 can be used to construct spherical 3-designs. We also conjecture that our lower bound gives the true value of $s(n, 3)$ (this has been verified for $n \leq 125$).

## 1. Introduction

We are interested in finding finite "well balanced" point sets on the surface of the unit $d$-sphere $S^d \subset R^{d+1}$. While it may be clear that vertices of regular polygons form such sets on the circle $S^1$, there is no natural way to generalize this for $d \geq 2$. Of the numerous possible criteria for measuring how "well balanced" our point set is (see e.g. [10]), one of the most useful and interesting one is that of the spherical design, as introduced in a monumental paper by Delsarte, Goethals, and Seidel in 1977 [11].

A *spherical t-design* on $S^d$ is a finite set of points $X \subset S^d$ for which the Chebyshev-type quadrature formula

$$\frac{1}{\sigma_d(S^d)} \int_{S^d} f(x) d\sigma_d(x) \approx \frac{1}{|X|} \sum_{x \in X} f(x)$$

is exact for all polynomials $f(x) = f(x_0, x_1, \ldots, x_d)$ of degree at most $t$ ($\sigma_d$ denotes the surface measure on $S^d$). In other words, $X$ is a spherical $t$-design of $S^d$, if for every polynomial $f(x)$ of degree $t$ or less, the average value of $f(x)$ over the whole sphere is equal to the arithmetic average of its values on the finite set $X$. General references on spherical designs include [11], [6], [5], and [22].

The existence of spherical designs for every $t$, $d$, and large enough $n = |X|$ was first proved by Seymour and Zaslavsky in 1984 [25], and general constructions were first given by the author in 1990 [3].

In [11], Delsarte, Goethals, and Seidel also proved that a spherical $t$-design on $S^d$ must have cardinality

$$n \geq N_d(t) = \binom{\lfloor t/2 \rfloor + d}{d} + \binom{\lfloor (t-1)/2 \rfloor + d}{d}.$$

A spherical $t$-design on $S^d$ with cardinality $N_d(t)$ is called *tight*. In 1980 Bannai and Damerell [7], [8] proved that tight spherical designs for $d \geq 2$ exist only for $t = 1, 2, 3, 4, 5, 7$ or $11$. All tight $t$-designs are known, except for $t = 4, 5$, and $7$. In particular, there is a unique tight spherical 11-design ($d = 23$ and $n = 196,560$).

Let $M_d(t)$ denote the minimum size of a spherical $t$-design on $S^d$, and let $M'_d(t)$ denote the smallest integer such that for every $n \geq M'_d(t)$, $t$-designs on $S^d$ exist on $n$ nodes. We have $N_d(t) \leq M_d(t) \leq M'_d(t)$. Values of $M_d(t)$ and $M'_d(t)$ are generally unknown when $d \geq 2$ and $t \geq 3$. For an upper bound on $M_d(t)$ and $M'_d(t)$ see [5].

The case $d = 1$ is completely settled; it is easy to see that vertices of a regular $n$-gon with $n \geq t + 1$ give a spherical $t$-design on the circle, hence $N_1(t) = M_1(t) + M'_1(t) = t + 1$. (Hong [19] proved in 1982 that these are the unique $t$-designs on $S^1$ when $t + 1 \leq n \leq 2t + 1$.)

Much work has been done for $d = 2$. It is well known that $N_2(t) = M_2(t)$ if and only if $t = 1$ (2 antipodal points), $t = 2$ (4 vertices of a regular tetrahedron), $t = 3$ (the regular octahedron), or $t = 5$ (the icosahedron). For $t = 4$ we have $N_2(4) = 9$, and there are designs of sizes $n = 12, 14$, and $n \geq 16$ [17]. Hardin and Sloane [17] also exhibit numerical evidence that a 4-design on $S^2$ does not exist for $n = 10, 11, 13$, and $15$; hence the conjectures $M_2(4) = 12$ and $M'_2(4) = 16$. Recent papers of Reznick [23] and Hardin and Sloane [18] give constructions for $t = 5$ (in which case $N_2(5) = M_2(5) = 12$) for $n = 12, 16, 18, 20$, and $n \geq 22$, and conjecture that this list is complete, hence that $M'_2(5) = 22$. In [18] Hardin and Sloane also provide numeric evidence for what they believe is a complete set of possible sizes for $t = 6, 7, 8, 9, 10, 11$, and $12$. Their work indicates that for these values of $t$, $M'_2(t) - M_2(t)$ varies greatly between 2 ($t = 12$) and 12 ($t = 7$).

Keeping $t$ constant and letting the dimension vary, we first note that $N_d(1) = M_d(1) = M'_d(1) = 2$ for every $d \geq 1$. Mimura [21] settled the case $t = 2$ in 1990: He proved that $M_d(2) = N_d(2) = d + 2$, and that $M'_d(2) = d + 2$ when $d$ is odd and $M'_d(2) = d + 4$ when $d$ is even. Much less has been known when $t \geq 3$. For $t = 3$ the author conjectured that 3-designs on $S^2$ do not exist on $n = 7$ or 9 points ($N_2(3) = 6$), and this was recently supported by a powerful computer search done by Hardin and Sloane [18]. In [17] Hardin and Sloane also present numerical evidence for values of $M_d(4)$ and $M'_d(4)$ for $d \leq 7$. If their conjectures are valid, then $M_d(4) = M'_d(4)$ for $d = 3, 4, 6$, and $7$, but $M'_d(4) - M'_d(4) = 12$ for $d = 5$.

The goal of this paper is to provide constructions for 3-designs on $S^d$ for all values of $n$ for which such designs exist of size $n$. Our results are summarized in the table below.

| $d$ | $NN_d(3) = M_d(3)$ | $n$ |
|---|---|---|
| 1 | 4 | $\geq 4$ |
| 2 | 6 | 6, 8, $\geq 10$ |
| 3 | 8 | 8, $\geq 10$ |
| 4 | 10 | 10, 12, $\geq 14$ |
| 5 | 12 | 12, $\geq 14$ |
| 6 | 14 | 14, 16, $\geq 18$ |
| 7 | 16 | 16, 18, $\geq 20$ |
| 8 | 18 | 18, 20, $\geq 22$ |
| 9 | 20 | 20, 22, $\geq 24$ |
| $\geq 5$ | $2d+2$ | $\geq 2d+2$ & even, $\geq 5(d+1)/2$ & odd |

We believe that our list above is complete. In particular, we conjecture that $M'_d(3) = \lfloor 5d/2 + 3 \rfloor_2$, where $d \neq 2$ or 4 and $\lfloor x \rfloor_2$ is the largest even integer not greater than $x$.

We will employ methods similar to those used in [1], [21], and [23]. We will also introduce and apply a concept from additive number theory generalizing the famous but not yet completely understood Sidon-sequences. A *Sidon-sequence*, as first studied by Sidon in 1993 [24], is a sequence of distinct integers $\{x_1, x_2, \ldots\}$ with the property that the sums $x_i + x_j$ are all distinct or, equivalently, that the equation $x_i + x_j - x_k - x_l = 0$ is satisfied only in the trivial case of $\{i,j\} = \{k,l\}$.

It follows from a 1941 paper of Erdös and Turán [14] (and was independently proved by Lindström in 1969 [20]) that in the interval $[1,n]$, a Sidon-sequence can have at most $n^{1/2} + n^{1/4} + 1$ elements. In 1944 Erdös [12] and Chowla [9] independently proved that a Sidon-sequence in $[1,n]$ with at least $n^1 - n^{5/16}$ elements can indeed be found. It is a \$1,000 Erdös problem to prove or disprove that the correct maximal cardinality differs from $\sqrt{n}$ by a constant. These and other results on Sidon-sets and related questions can be found in Erdös's and Freud's excellent survey [13], as well as in [15] and [16].

In this paper we are interested in the following generalization. Let $S$ be a set of integers, and suppose that the congruence $\varepsilon_1 x_1 + \varepsilon_2 x_2 + \cdots + \varepsilon_t x_t \equiv 0 \mod n$, where $\varepsilon_i = 0, \pm 1$ and $x_i \in S$ for $i = 1, 2, \ldots, t$, only holds in the trivial case, that is when $\varepsilon_i = 0$ for all $i = 1, 2, \ldots, t$ or when the same $x_i$ appears with both a coefficient of 1 and of $-1$. We here call such sets *Sidon-type sets of strength t*, and denote their maximum cardinality (they clearly must be finite) by $s(n,t)$. It is obvious that $s(n,1) = n-1$, and it is also easy to see that $s(n,2) = \lfloor (n-1)/2 \rfloor$. Here we find the following lower bound for $s(n,3)$: (i) $s(n,3) \geq \lfloor n/4 \rfloor$ is $n$ is even; (ii) $s(n,3) \geq \lfloor (n+1)/6 \rfloor$ if $n$ is odd and has no divisors congruent to 5 mod 6; and (iii) $s(n,3) \geq \frac{(p+1)n}{6p}$ if $n$ is odd and $p$ is its smallest divisor which is congruent to 5 mod 6. We show how Sidon-type sets of strength 3 can be used to construct spherical 3-designs. We also conjecture that our lower bound gives the true value of $s(n,3)$ (this has been verified for $n \leq 125$), which in part supports our conjecture for $M'_d(3)$ above. Note also that a Sidon-type set of strength 4 forms a Sidon-sequence in $[1, n]$, hence $s(n, 4) \leq n^{1/2} + n^{1/4} + 1$.

## 2. Harmonic Polynomials

To construct spherical designs, we will use the following equivalent definition, cf. [11]:

*A finite subset $X$ of $S^d$ is a spherical t-design, if and only if*

$$\sum_{x \in X} f(x) = 0$$

*for all homogeneous harmonic polynomials $f(x_0, x_1, \ldots, x_d)$ with $1 \leq \deg f \leq t$.*

A polynomial $f(x_0, x_1, \ldots, x_d)$ is called *harmonic* if it satisfies Laplace's equation $\Delta f = 0$. The set of homogeneous harmonic polynomials of degree $s$ forms a vector space $Harm_{d+1}(s)$, with

$$\dim Harm_{d+1}(s) = \binom{s+d}{d} - \binom{s+d-2}{d}.$$

In particular, for $s \leq 3$, we see that $\Phi_s$ forms a basis for $Harm_{d+1}(s)$ where

$\Phi_1 = \{x_i | 0 \leq i \leq d\}$,

$\Phi_2 = \{x_i x_j | 0 \leq i < j \leq d\} \cup \{x_i^2 - x_{i+1}^2 | 0 \leq i \leq d-1\}$, and

$\Phi_3 = \{x_i x_j x_k | 0 \leq i < j < k \leq d\} \cup \{x_i^3 - 3x_i x_j^2 | 0 \leq i \neq j \leq d\}$.

We associate matrices with spherical designs in the following way. For a set $X = \{u_k = (u_{ok}, u_{1k}, \ldots, u_{dk}) \in R^{d+1} | 1 \leq k \leq n\}$ we consider the $(d+1) \times n$ matrix $U$ with column vectors $u_1, u_2, \ldots, u_n$.

For a polynomial $f(x_0, x_1, \ldots, x_d)$ we define

$$f(U) = \sum_{k=1}^{n} f(u_k).$$

With these notations, $X$ is a spherical $t$-design, if and only if

(∗)
$$\sum_{i=0}^{d} u_{ik}^2 = 1 \quad \text{for } 1 \leq k \leq n, \quad \text{and}$$

$$f(U) = 0 \quad \text{for every polynomial } f \in \bigcup_{s=1}^{t} \Phi_s.$$

## 3. Antipodal Designs

It is well known and most obvious that vertices of the generalized regular octahedra form (tight) 3-designs on $S^d$:

**Construction 3.1.** *The matrix $(I\ -I)$ provides a spherical 3-design on $S^d$ of size $2d + 2$. Here $I$ is the $d + 1$ by $d + 1$ identity matrix.*

More generally, antipodal point sets on $S^d$ (sets where $x \in S^d$ implies $-x \in S^d$) can be used to construct spherical 3-designs. Equations (*) show that if $t$ is even and $A$ is the matrix of a $t$-design on $S^d$, then (the set of column vectors of) the matrix $(A\ -A)$ provides a $(t+1)$-design on $S^d$. Since 2-designs on $S^d$ exist for sizes $n \geq d + 2$ when $d$ is odd and for $n = d + 2, n \geq d + 4$ when $d$ is even [21], we immediately have

**Proposition 3.2.** *Let $n$ be an even integer such that $n \geq 2d + 4$, except for $n = 2d + 6$ when $d$ is even. Then a spherical 3-design on $S^d$ of size $n$ exists.*

Primarily with the cases of even $d$ in mind, we provide the following

**Construction 3.3.** *Suppose that $A$ is the matrix of a 2-design on $S^{d-1}$ of size $n_1$, $J$ is the 1 by $n_1$ matrix of all 1's, $\alpha = \sqrt{d/(d+1)}$, and $\delta = \sqrt{1/(d+1)}$. Then $M = \begin{pmatrix} \alpha A & -\alpha A \\ \delta J & -\delta J \end{pmatrix}$ is a 3-design of size $2n_1$ on $S^d$.*

*Proof.* For $A = (u_{ik})_{0 \leq i \leq d-1, 1 \leq k \leq n_1}$ we have

$$\sum_{i=0}^{d-1} u_{ik}^2 = 1 \quad \text{for } 1 \leq k \leq n_1, \quad \text{and}$$

$$\sum_{k=0}^{n_1} u_{ik}^2 - u_{i+1,k}^2 = 0 \quad \text{for } 0 \leq i \leq d - 2.$$

Therefore, $M$ satisfies (*) for $t = 3$ if and only if the equations

$$\alpha^2 + \delta^2 = 1 \quad \text{and} \quad \alpha^2 \frac{2n_1}{d} - \delta^2 \cdot 2n_1 = 0$$

hold. These two equations are equivalent to

$$\alpha^2 = \frac{d}{d+1} \quad \text{and } \delta^2 = \frac{1}{d+1}. \qquad \square$$

As a corollary, we get

**Proposition 3.4.** *Let $n$ be an even integer such that $n \geq 2d + 2$, except for $n \geq 2d + 4$ when $d$ is odd. Then a spherical 3-design on $S^d$ of size $n$ exists.*

## 4. Sidon-Type Sets

For other constructions of spherical 3-designs, we will use what we call Sidon-type sets of strength 3.

Let $R$ be a ring with identity, $S$ a subset of $R$, and $t$ a positive integer. We say that $S$ is a *Sidon-type set of strength $t$* in $R$ if no non-trivial-trivial sum of the

form

$$\varepsilon_1 x_1 + \varepsilon_2 x_2 + \cdots + \varepsilon_t x_t,$$

where $\varepsilon_1, \varepsilon_2, \ldots, \varepsilon_t = 0, \pm 1$ and $x_1, x_2, \ldots, x_t$ are (not necessarily distinct) elements of $S$, equals 0. We call such a sum non-trivial if no $x_i$ appears in it with both a coefficient of 1 and $-1$, and if at least one $\varepsilon_i$ is non-zero ($i = 1, 2, \ldots, t$).

Here we are only interested in Sidon-type sets in $Z_n$, and we think of these sets as integer subsets of the interval $[1, n]$. The cardinality of a largest Sidon-type set of strength $t$ in $Z_n$ will be denoted by $s(n, t)$. It is obvious that $s(n, 1) = n - 1$ (take all integers from 1 to $n - 1$), and it is easy to see that $s(n, 2) = \lfloor (n - 1)/2 \rfloor$ ($S$ cannot contain both $x$ and $n - x$, but it can consist of all integers up to $\lfloor (n - 1)/2 \rfloor$). For $t = 3$ we give a constructive proof for the following.

**Theorem 4.1.**

(i)  $s(n, 3) \geq \lfloor n/4 \rfloor$ if $n$ is even;
(ii)  $s(n, 3) \geq \lfloor (n + 1)/6 \rfloor$ if $n$ is odd and has no divisors congruent to 5 mod 6; and
(iii)  $s(n, 3) \geq \frac{(p+1)n}{6p}$ if $n$ is odd and $p$ is its smallest divisor which is congruent to 5 mod 6.

*Proof.* We can always take all the odd integers up to (but not including) $n/3$, proving (ii). When $n$ is even, we can take all the odd integers up to (but not including) $n/2$, which proves (i).

Now suppose that $n$ is odd and that there is a non-negative integer $q$ such that $p = 6q + 5$ divides $n$. We define

$$S = \{ap + 2b + 1 \mid a = 0, 1, \ldots, n/p - 1, \quad b = 0, 1, \ldots, q\}.$$

We see that $S$ has cardinality $(p + 1)n/(6p)$. To verify that $S$ is a Sidon-type set of strength 3, suppose that $n$ divides

$$x = \varepsilon_1 x_1 + \varepsilon_2 x_2 + \varepsilon_3 x_3$$
$$= (\varepsilon_1 a_1 + \varepsilon_2 a_2 + \varepsilon_3 x_3)p + 2(\varepsilon_1 b_1 + \varepsilon_2 b_2 + \varepsilon_3 b_3) + \varepsilon_1 + \varepsilon_2 + \varepsilon_3.$$

This implies that

$$y = 2(\varepsilon_1 b_1 + \varepsilon_2 b_2 + \varepsilon_3 x_3) + \varepsilon_1 + \varepsilon_2 + \varepsilon_3$$

is divisible by $p$, but since $|y| \leq 6q + 3$ and $p = 6q + 5$, this can only happen if $y = 0$. Since 0 is an even number, either all $\varepsilon$'s are equal to 0 (a trivial sum), or exactly one $\varepsilon$, say $\varepsilon_1$, is 0. In the latter case, since $b_2, b_3 \geq 0$, we must have $\varepsilon_2 = -\varepsilon_3$, which implies that $b_2 = b_3$. In this case we also get

$$|X| = |(\varepsilon_1 a_1 + \varepsilon_2 a_2 + \varepsilon_3 a_3)p| = |a_2 - a_3|p \leq \left(\frac{n}{p} - 1\right)p < n,$$

so $n$ can only divide $x$ if $x = 0$. But then $a_2 = a_3$, hence $x_2 = x_3$, and we again have the trivial sum. □

We performed a computer search for values of $s(n, 3)$ for $n \leq 125$, and found that in all cases the true value agreed with the lower bound found in Theorem 4.1. Therefore we state

**Conjecture 4.2.** *Theorem* 4.1 *gives the exact value of $s(n, 3)$. In particular,*

(i) $s(n, 3) \leq n/4$, *with equality if and only if $n$ is divisible by* 4; *and*
(ii) *if $n$ is odd, then $s(n, 3) \leq n/5$, with equality if and only if $n$ is divisible by* 5.

## 5. Regular 3-Designs

It is well known and is easy to check that the vertices of a regular $n$-gon where $n \geq t + 1$ form a $t$-design on $S^1$. In this section we investigate a generalization of this for dimensions $d \geq 1$, where $d$ is odd.

For positive integers $m$ and $n$ we define the vectors $s(m)$ and $c(m)$ in $R^n$ to be

$$s(m) = \left(\sin\left(\frac{2\pi}{n}m\right), \sin\left(\frac{2\pi}{n}2m\right), \ldots, \sin\left(\frac{2\pi}{n}nm\right)\right) \text{ and}$$

$$c(m) = \left(\cos\left(\frac{2\pi}{n}m\right), \cos\left(\frac{2\pi}{n}2m\right), \ldots, \cos\left(\frac{2\pi}{n}nm\right)\right).$$

Now let $e > 0$ and $m_1, m_2, \ldots, m_e$ be integers, and set $S = \{m_1, m_2, \ldots, m_e\}$. We define the matrix $A(S)$ to be the $2(e) \times n$ matrix with row vectors $s(m_1), c(m_1), s(m_2), c(m_2), \ldots, s(m_e), c(m_e)$.

**Lemma 5.1.** *Let $e$ be a positive integer, $s = 1, 2$, or 3, and suppose that $S = \{m_1, m_2, \ldots, m_e\}$ is a Sidon-type set of strength $s$. We define the matrix $A = A(S)$ as above. If $f : R^{2e} \to R$ is a polynomial such that $f \in \Phi_s$, then $f(A) = 0$.*

*Proof.* The statement is well known for $s = 1$: setting $z_j = \cos(\frac{2\pi}{n}m_j) + i\sin(\frac{2\pi}{n}m_j)$, we see that, since $z_j \neq 1$ when $m_j$ is not divisible by $n$, we have $\sum_{k=1}^{n} z_j^k = 0$ for every $j = 1, 2, \ldots, e$.

Our claims for the cases of $s = 2$ when $f$ is square-free and of $s = 3$ (and $f$ any homogeneous cubic polynomial) follow from the $s = 1$ case after repeated use of the trigonometric identities

$$\sin x \sin y = \tfrac{1}{2}[\cos(x - y) - \cos(x + y)],$$
$$\sin x \cos y = \tfrac{1}{2}[\sin(x + y) + \sin(x - y)], \text{ and}$$
$$\cos x \cos y = \tfrac{1}{2}[\cos(x + y) + \cos(x - y)].$$

Finally, when $f = x_i^2 - x_{i+1}^2$, $i = 0, 1, \ldots, d - 1$, we get

$$f(A) = \frac{1}{2}\sum_{k=1}^{n} \cos 0 - \frac{1}{2}\sum_{k=1}^{n} \cos 0 = 0. \qquad \square$$

We have the following corollary.

**Construction 5.2.** *Suppose that d is an odd positive integer, $e = (d + 1)/2$, $s = 1, 2$, or 3, and that $S = \{m_1, m_2, \ldots, m_e\}$ is a Sidon-type set of strength s. Then the n column vectors of the matrix $M(S) = \sqrt{2/(d+1)} \cdot A(S)$ form a spherical s-design on $S^d$.*

The *s*-designs constructed with Construction 5.2 will be called *regular s-designs*. We note that Lemma 5.1 and Construction 5.2 are false for strengths $s \geq 4$ (see [1]). As Construction 5.2 provides spherical 3-designs on $S^d$ with $s(n, 3) \geq (d + 1)/2$, we have

**Proposition 5.3.** *Suppose that d is an odd positive integer. Regular 3-designs of size n on $S^d$ exist when*

(i)  *n is even and $n \geq 2d + 2$;*
(ii) *n is odd and $n \geq 3d + 2$; and*
(iii) *n is odd and $n \geq \frac{p}{p+1}(3d + 3)$, where p is a divisor of n which is congruent to 5 mod 6.*

*In particular*, there are regular 3-*designs of size n on $S^d$ when n is an odd integer which is divisible by 5 and $n \geq 5(d + 1)/2$.*

Conjecture 4.2 implies that Proposition 5.3 characterizes all values of *n* for which regular 3-designs exist on $S^d$. In particular, we believe that no regular 3-design exists for odd values of *n* with $n < 5(d + 1)/2$ (this has been verified for $d \leq 49$).

## 6. Other Spherical 3-Designs

We have just seen constructions for 3-designs on $S^d$ for all odd values of *n* when $n \geq 5(d + 1)/2$, *d* is odd, and *n* is divisible by 5. In this section we will construct 3-designs on $S^d$ of size *n* for every odd value of *n* with $n \geq \max\{5(d + 1)/2, 2d + 7\}$.

**Construction 6.1.** *Let $d_1$ and $d_2$ be positive even integers with $d = d_1 + d_2 - 1$, let $n_1$ and $n_2$ be positive integers with $n = n_1 + n_2$, and suppose that $d_1 n_1 \geq d_2 n_2$. Suppose further that A is (the matrix of, see section 2) a regular 3-design of size $n_1$ on $S^d$, and that C is a regular 3-design of size $n_2$ on $S^{d_1-1}$. Then the n column vectors of the matrix $M = \begin{pmatrix} \alpha A_1 & C \\ \beta A_2 & 0 \end{pmatrix}$, where $A = \begin{pmatrix} A_1 \\ A_2 \end{pmatrix}$ and $\alpha$ and $\beta$ are defined below, form a 3-design on $S^d$.*

*Proof.* $M$ satisfies (*) if and only if the equations

$$\alpha^2 \frac{d_1}{d_1 + d_2} + \beta^2 \frac{d_2}{d_1 + d_2} = 1 \text{ and}$$

$$\frac{n_2}{d_1} + \alpha^2 \frac{n_1}{d_1 + d_2} - \beta^2 \frac{n_1}{d_1 + d_2} = 0$$

hold (see proof of Construction 3.3). The two equations are equivalent to

$$\alpha^2 = 1 - \frac{d_2 n_2}{d_1 n_1} \text{ and } \beta^2 = 1 + \frac{n_2}{n_1}.$$

We see that $\alpha$ is real iff $d_1 n_1 \geq d_2 n_2$, as was assumed. □

A corollary is the following.

**Proposition 6.2.** *Let $d$ and $n$ be odd integers such that $n \geq \max\{5(d+1)/2, 2d+7\}$. Then there are 3-designs of size $n$ on $S^d$.*

*Proof.* The construction is given in 6.1 when taking $d_1 = 2$, $d_2 = d - 1$, $n_1 = n - 5$, and $n_2 = 5$. The necessary inequalities all hold: $d_1 n_1 \geq d_2 n_2$, $n_1 \geq 2d + 2$ ($n_1$ is even, see Proposition 5.3 (i)), and $n_2 \geq 4$ (C is a 3-design on the circle, see Proposition 5.3 (ii)). □

We now turn to the case when $d$ is even.

**Construction 6.3.** *Let $d_1$ and $d_2$ be positive even integers with $d = d_1 + d_2$, let $n_1$ and $n_2$ be positive integers with $n = 2n_1 + n_2$, and suppose that $2d_1 n_1 \geq (d_2 + 1)n_2$. Suppose further that $A$ is a regular 2-design of size $n_1$ on $S^{d-1}$, and that $C$ is a regular 3-design of size $n_2$ on $S^{d_1-1}$. Then the $n$ column vectors of the matrix $M = \begin{pmatrix} \alpha A_1 & -\alpha A_1 & C \\ \beta A_2 & -\beta A_2 & 0 \\ \delta J & -\delta J & 0 \end{pmatrix}$, where $A = \binom{A_1}{A_2}$, $J$ is the 1 by $n_1$ matrix of all 1's, and $\alpha, \beta,$ and $\delta$ are defined below, form a 3-design on $S^d$.*

*Proof.* $M$ satisfies (*) if and only if the equations

$$\alpha^2 \frac{d_1}{d_1 + d_2} + \beta^2 \frac{d_2}{d_1 + d_2} + \delta^2 = 1,$$

$$\frac{n_2}{d_1} + \alpha^2 \frac{2n_1}{d_1 + d_2} - \beta^2 \frac{2n_1}{d_1 + d_2} = 0,$$

$$\frac{n_2}{d_1} + \alpha^2 \frac{2n_1}{d_1 + d_2} - \delta^2 2n_1 = 0, \text{ and}$$

$$\beta^2 \frac{2n_1}{d_1 + d_2} - \delta^2 2n_1 = 0$$

hold (see again the proof of Construction 3.3). The four equations are equivalent to the following three:

$$\alpha^2 = \frac{d_1 + d_2}{d_1 + d_2 + 1} \left(1 - \frac{(d_2 + 1)n_2}{2d_1 n_1}\right),$$

$$\beta^2 + \frac{d_1 + d_2}{d_1 + d_2 + 1} \left(1 + \frac{n_2}{2n_1}\right), \text{ and}$$

$$\delta^2 = \frac{1}{d_1 + d_2 + 1} \left(1 + \frac{n_2}{2n_1}\right).$$

We see that $\alpha$ is real iff $2d_1 n_1 \geq (d_2 + 1)n_2$, as was assumed. □

This gives us the following corollary.

**Proposition 6.4.** *Let $d$ be even and $n$ odd, such that $n \geq \max\{5(d+1)/2, 2d+7\}$. Then there are 3-designs of size $n$ on $S^d$.*

*Proof.* The construction is given in 6.3 when taking $d_1 = 2$, $d_2 = d - 2$, $n_1 = (n-5)/2$, and $n_2 = 5$. The necessary inequalities all hold: $2d_1 n_1 \geq (d_2 + 1)n_2$, $n_1 \geq d + 1$ (Construction 5.2 provides regular 2-designs on $S^{d-1}$ when $d - 1$ is odd and $s(n_1, 2) = \lfloor (n_1 - 1)/2 \rfloor \geq d/2$), and finall $n_2 \geq 4$ (see Proposition 5.3(ii)). □

## 7. Summary of Results

We summarize our results in the following theorem.

**Theorem 7.1.** *No spherical 3-design exists on $S^d$ of size $n < 2d + 2$. Spherical 3-designs on $S^d$ exist when*

(i) *$n$ is even and $n \geq 2(d+1)$;*
(ii) *$n$ is odd and $n \geq 5(d+1)/2$, except for $d = 2$ and $n = 9$, $d = 4$ and $n = 13$.*

*Proof.* The first statement is a reiteration of $n \geq N_d(3) = 2d + 2$ (see [11] for more). The cases when $n$ is even are either from Proposition 3.2 or Proposition 3.4. The cases of $d = 1$ and $d = 3$ of (ii) are stated in Proposition 5.3 (ii), and the case of $d = 5$ and $n = 15$ is a special case of Proposition 5.3 (iii). All other cases in (ii) follow from Propositions 6.2 and 6.4. □

According to Theorem 7.1, the number of different values of $n$ for which the problem is open is $\lfloor (d+2)/4 \rfloor$ when $d \neq 2$ or 4. Only one case is unsettled for $d = 3$ and $d = 5$, and at most two cases are open when $d \leq 9$. We state our

**Conjecture 7.2.** *Theorem 7.1 gives a complete list of all possible sizes $n$ for which spherical 3-design on $S^d$ exist. In particular, $M'_d(3) = \lfloor 5d/2 + 3 \rfloor_2$, where $d \neq 2$ or 4 and $\lfloor x \rfloor_2$ is the largest even integer not greater than $x$, $M'_2(3) = 10$, and $M'_4(3) = 14$.*

Conjecture 7.2 is supported by our previously stated belief that no regular 3-designs exist with $n < 5(d+1)/2$ when $n$ is odd (which has been verified for $d \leq 49$); and was also numerically demonstrated for $d = 2$ by Hardin and Sloane [18].

**Acknowledgments.** I would like to thank Charles A. Ross for his computer justification of Conjecture 4.2, and Róbert Freud and Bruce Reznick for helpful comments.

*Note added in proof.* An upcoming paper of Boyvalenkov, Danev, and Nikova contains new nonexistence results, such as the nonexistence of a 7 point 3-design on $S^2$.